\documentclass[11pt,twoside]{article}

\usepackage{amsbsy,amsfonts,amsmath,amssymb,eucal,mathrsfs,amscd}
\usepackage{diagrams}

\def\v#1{\par \vspace{#1mm} \par}

\def\qed{\hfill $\square$ \vspace{5mm}}

\renewcommand{\tilde}{\widetilde}
\renewcommand{\hat}{\widehat}
\renewcommand{\bar}{\overline}

\newcommand{\Z}{\mathbb Z}

\newcommand{\C}{\mathbb{C}}

\newcommand{\0}{\mathcal O}
\newcommand{\g}{\mathfrak g}

\renewcommand{\u}{\mathfrak u}
\newcommand{\n}{\mathfrak n}
\newcommand{\e}{\mathfrak e}
\newcommand{\plie}{\mathfrak p}
\newcommand{\q}{\mathfrak q}

\newcommand{\p}{{\mathbb P}}

\renewcommand{\S}{{\mathbb S}}

\newcommand{\co}{{\cal O}}

\renewcommand{\iff}{if and only if }

\newcommand{\tr}{{}^t}

\newcommand{\scal}[1]{\langle #1 \rangle}
\newcommand{\im}{\mathtt{Im}}
\newcommand{\lpara}{
\ \vspace{3pt}

\noindent}

\newcommand{\liste}{\

\begin{itemize}}
\newcommand{\codim}{\mbox{codim}}

\newcommand{\fonctionrat}[5]{
\begin{array}[t]{rrcll}
#1 & : & #2 & \dasharrow & #3 \\
   &   & #4 & \mapsto     & #5
\end{array}  }

\newtheorem{definition}{Definition}
\newenvironment{defi}{\begin{definition} \rm}{\end{definition}}

\newtheorem{prop}{Proposition}

\newtheorem{lemm}{Lemma}
\newtheorem{coro}{Corollary}
\newtheorem{theo}{Theorem}
\newtheorem{notation}[definition]{Notation}

\newtheorem{construction}[definition]{Construction}

\newtheorem{rem}{Remark}

\newtheorem{exam}{Example}

\newtheorem{examples}[definition]{Examples}

\newtheorem{nothing}[definition]{$\!\!$}

\newenvironment{proo}{{\flushleft \bf Proof :}}{\hfill $\square$ \vspace{5mm}}

\newtheorem{definition*}{Definition}[section]
\newenvironment{defi*}{\begin{definition*} \rm}{\end{definition*}}
\newtheorem{definitions*}[definition*]{Definitions}
\newenvironment{defis*}{\begin{definitions*} \rm}{\end{definitions*}}
\newtheorem{prop*}[definition*]{Proposition}
\newtheorem{lemm*}[definition*]{Lemma}
\newtheorem{coro*}[definition*]{Corollary}
\newtheorem{theo*}[definition*]{Theorem}
\newtheorem{remark*}[definition*]{Remark}
\newenvironment{rema*}{\begin{remark*} \rm}{\end{remark*}}
\newtheorem{remarks*}[definition*]{Remarks}
\newenvironment{remas*}{\begin{remarks*} \rm}{\end{remarks*}}
\newtheorem{example*}[definition*]{Example}
\newenvironment{exam*}{\begin{example*} \rm}{\end{example*}}
\newtheorem{examples*}[definition*]{Examples}
\newenvironment{exams*}{\begin{examples*} \rm}{\end{examples*}}
\newtheorem{nothing*}[definition*]{$\!\!$}
\newenvironment{noth*}{\begin{nothing*} \rm}{\end{nothing*}}

\newfont{\joli}{cmfi10 scaled \magstep2}

\DeclareMathAlphabet{\mathpzc}{OT1}{pzc}{m}{it}

   \def\cN{{\cal N}}
\def\cO{{\cal O}}
 \def\cQ{{\cal Q}}

\begin{document}

\begin{center}
{\bf \Large On stratified Mukai flops}  \v{10}{\bf Pierre-Emmanuel
Chaput and Baohua Fu} \v{2}
\end{center}

\v{7}

{\def\thefootnote{\relax}
\footnote{ \hspace{-6.8mm}
}
}
\section{Introduction}

In  recent studies (see \cite{Na2}, \cite{fu}) of birational
geometry of symplectic resolutions of nilpotent orbit closures,
three types of flops (which will be called stratified Mukai flops
of type $A, D, E_6$) are shown to be fundamental, in the sense
that others can be decomposed into a sequence of these flops.
Stratified Mukai flops of type $A$ are given by cotangent bundles
of dual Grassmanians $ T^*G(k, V) \dasharrow T^*G(k, V^*)$, which
has been previously studied by Markman(\cite{markman}). The type
$D$ stratified Mukai flop is given by the cotangent bundles of the
two connected components of the orthogonal Grassmannian
$G_{iso}(k, 2k)$, where $k$ is an odd integer. There are two
stratified Mukai flops of type $E_6$, corresponding to the pairs
of roots $(\alpha_1, \alpha_6)$ (type $E_{6, I}$) and $(\alpha_3,
\alpha_5)$(type $E_{6,II}$).

For a stratified Mukai flop $\mu: T^*(G/P) \dasharrow T^*(G/Q)$,
an important problem is to construct explicitly a $G$-equivariant equivalence
between the derived categories of $T^*(G/P)$ and $T^*(G/Q)$ which
is compatible with the flop $\mu$. It is tempting to believe that
a Fourier-Mukai transform with kernel induced by a suitable
resolution of this flop could yield such an equivalence.

In  \cite{chaput} and \cite{fu}, we showed that the graph of the
stratified Mukai flop is in fact smooth, thus it gives a
resolution of the flop. Although very natural, this resolution
does not lead to easy computations in Chow groups  nor in the
derived categories.
Our first result (Theorem \ref{resolution},
which works also for type $D$ Mukai flops) of this note is to
resolve the stratified Mukai flop of type $E_{6,I}$ by blow-ups
along smooth centers. A similar resolution has already been constructed
by Markman \cite{markman} for stratified Mukai flops of type $A$.
As we will see, our resolution is constructed in a very similar
way, blowing up successively the different $G$-orbits, from the
smallest to the biggest.

Like the usual Mukai flops, this resolution does not induce an
isomorphism between the Chow groups (cf. Example \ref{ex1}, Example \ref{ex2}).
To find out a natural
isomorphism between the Chow groups, we will first construct a
deformation of this flop, then we will show that the deformed flop
can be resolved in a similar way. Using a beautiful idea of
\cite{LLW}, we prove that
 this leads to a natural
isomorphism between Chow groups and their Chow motifs. This shows
that the stratified Mukai flop of type $E_{6, I}$ is very similar to
the usual Mukai flop. Unfortunately, we do not know if this
functor also gives an equivalence between derived categories.

 \vspace{0.3 cm}

{\em Acknowledgements:} We are grateful to E. Markman for his
clear explanations on  \cite{markman}.
  The second named author would like to thank Y. Namikawa for helpful
   discussions and the
Max Planck Institute (Bonn) for the hospitality.
\section{The wonderful resolution}
\label{section_wonderful}

Let $G$ be a simple group of type $A_n,D_{2n+1}$ or $E_6$.
Let $(P,Q)$ be the maximal parabolic subgroups corresponding respectively
to the roots
$(\alpha_i,\alpha_{n+1-i})$ with $2i<n+1$,
$(\alpha_{2n},\alpha_{2n+1})$ and $(\alpha_1,\alpha_6)$.

Set $X_P = G/P$ and  $T_P = T^* (G/P)$.  We denote by $\pi_P:T_P
\rightarrow X_P$ the natural projection. We know that the order on
the set of $G$-orbits in $T_P$ given by $\cO \leq \cO'$ if $\co
\subset \overline{\cO'}$ is a total order, so we choose a
labelling $T_P^i,0 \leq i\leq r$ of the $G$-orbits in $T_P$ such
that $T_P^i \subset \overline {T_P^j}$ \iff $i \leq j$.  Note that
$r=i$ in type $A_n$ with $P,Q$ corresponding to
$\alpha_i,\alpha_{n+1-i}$, and that $r$ equals $n$ in type
$D_{2n+1}$ and $2$ in type $E_6$. Similar notations will be used
for the parabolic subgroup $Q$. It is known that there exists a
birational map $\mu: T_P \dasharrow T_Q$, which will be called a
{\em stratified Mukai flop}. Note that $\mu$ is an isomorphism
between the open orbits $T_P^r \to T_Q^r$.

We consider successive blow-ups defined recursively as follows.
Let $Bl_0(T_P) := T_P$ and $Bl_0(\overline{T_P^i}) =
\overline{T_P^i}$. For $1 \leq k \leq r$ an integer, we let
$Bl_k(T_P)$ be the blow-up of $Bl_{k-1}(T_P)$ along
$Bl_{k-1}(\overline{T_P^{k-1}})$,  $Bl_k(\overline{T_P^{k-1}})$
the exceptional divisor of this blow-up and
 $Bl_k(\overline{T_P^i})\ (i \not = k-1)$
the proper transform of $Bl_{k-1}(\overline{T_P^i})$ under this
blow-up. We set $\tilde T_P := Bl_r(T_P)$ and $\tilde T_P^i :=
Bl_r(\overline{T_P^i})$. Note that the natural blow down map
restricts to an isomorphism over the open $G$-orbit. We also
consider the analogous construction for $Q$.

\begin{theo}            \label{resolution}
Under the above notations, we have

(i) The varieties $Bl_{k}(\overline{T_P^{k}})$, $Bl_k(T_P)$ are smooth
for any $k$. The divisors
$\tilde T_P^i (i=0, \cdots, r-1)$ are smooth irreducible
with normal crossing inside $\tilde T_P$ and
$Pic(\tilde T_P) \simeq Pic(G/P) \oplus_{i=0}^{r-1} \Z [\tilde T_P^i]  $;

(ii) Two points $x,y \in \tilde T_P$ belong to the same $G$-orbit \iff
$\{i \mid  x \in \tilde T_P^i \} = \{ j \mid y \in \tilde T_P^j \}$;

(iii) There is a unique $G$-equivariant isomorphism
$\tilde \mu : \tilde T_P \simeq \tilde T_Q$ such that
$\tilde \mu$ identifies with
$\mu$ over the open $G$-orbit.
\end{theo}

\begin{rem}
1). Claims $(i)$ and $(ii)$ are some properties of wonderful
compactifications of symmetric spaces (\cite{deconcini}). This is
why we called this resolution wonderful.

2). We do not know if a similar result holds for stratified Mukai
flop of type $E_{6,II}$. In this case, the inclusion relationship
of orbit closures in $T^*X_P$ is no longer linear.
\end{rem}

We will need the following lemma in the proof of the theorem.

\begin{lemm} \label{lift}
Let $Y$ be a variety and assume there are morphisms $t_P:Y
\rightarrow T_P,x_Q:Y \rightarrow X_Q$ such that $t_P^{-1}(T^r_P)$
is dense in $Y$ and $\forall y \in t_P^{-1}(T^r_P), \pi_Q \circ
\mu \circ t_P(y)  = x_Q(y)$. Then there is a unique morphism $t_Q
: Y \rightarrow T_Q$ which lifts $x_Q$ via $\pi_Q$ and  such that
$\forall y \in t_P^{-1}(T^r_P), \mu \circ t_P (y) = t_Q(y)$.
\end{lemm}
\begin{proo}
The uniqueness of $t_Q$ is clear. Now we prove the existence
following \cite{chaput}.  A pair $(p,q) \in X_P \times X_Q$ is
called {\em incident} if $Stab(p) \cap Stab(q)$ is a parabolic
subgroup of $G$. For $p \in X_P$, let $C_p \subset X_Q$ denote the
set of points incident to $p$, and define similarly $C_q$, for $q
\in X_Q$. For an incident pair $(p,q)$, recall that
there is a well-defined linear isomorphism
$\mu(p,q):T_p X_P/T_p C_q \rightarrow T_q X_Q / T_q C_p$
\cite[theorem 4.1]{chaput}. Let ${(T_q C_p)}^\bot$ be
the subspace of $T^*_qX_Q$ consisting of co-vectors vanishing
identically on $T_q C_p$. Similarly, we get ${(T_p C_q)}^\bot$.
 The transposed map $\tr
\mu(p,q) : {(T_q C_p)}^\bot \rightarrow {(T_p C_q)}^\bot$ is then
again an isomorphism. Furthermore by \cite[theorem 4.1]{chaput},
for any  $t \in {(T_p C_q)}^\bot \cap T^r_P$, we have $\mu(t) = (q
, \tr \mu(p,q)^{-1}(t) )$.

Note that for any $y \in Y$, the pair $p:=\pi_P \circ t_P(y)$ and
$q:=x_Q(y)$ are incident. In fact, this is true if $y \in
t_P^{-1}(T^r_P)$(see for example \cite[theorem 3.3]{chaput}). Now
the density of $t_P^{-1}(T^r_P)$ implies that this holds for all
$y \in Y$. We claim that for any $y \in Y$, $t_P(y) \in {(T_p
C_q)}^\bot.$ In fact, if $y \in t_P^{-1}(T^r_P),$ then $ t_P(y)
\in T_p^*X_P \cap T_q^*X_Q =  \plie^\bot \cap \q^\bot\simeq
(T_pC_q)^\bot,$ where $\plie$ and $\q$ are  the Lie subalgebras of
$Stab(p)$ and $Stab(q)$. Again the claim follows from the density
of $t_P^{-1}(T^r_P)$. Therefore, to conclude the proof of the
lemma, it is enough to set $x_P(y) := \pi_P \circ t_P(y)$ and
define
$$t_Q(y) := ( x_Q(y), \tr \mu(x_P(y),x_Q(y))^{-1}(t_P(y)) ). $$
\end{proo}

\noindent {\bf Proof of the theorem : } The uniqueness  of $\tilde
\mu$ follows from the fact that $\tilde T_P^r$ is dense in $\tilde
T_P$, so let us prove the existence of $\tilde \mu$. Our proof is
in three steps~: first we give a description of $\tilde T_P$ in
terms of complete collineations; from this description we then
deduce the existence of a morphism $\tilde T_P \rightarrow \tilde
T_Q$. Finally, we show that this morphism and the analogous
morphism $\tilde T_Q \rightarrow \tilde T_P$ are inverses one to
the other. This method is very similar to that of Markman
(\cite{markman}), where he proved the theorem for stratified Mukai
flops of type A. We will give a complete proof in the case of
$E_6$, and explain how this adapts to type $D$.

So consider first the case $G = E_6$. Then a Levi factor $L$ of
$P$ is a semi-direct product of $\C^*$ and $Spin_{10}$, and if $x$
denotes the point in $G/P$ with stabilizer $P$, the representation
$T^*_xX_P$ of $L$ is a (16-dimensional) spin representation
\cite[theorem 2.10]{faulkner}. We will denote this
$L$-representation by $\S_+$.

\smallskip

Let $W$ be the natural
10-dimensional
$L$-representation and $\cQ \subset \p W$ the corresponding smooth
quadric. The variety  $\p \S_+$ has two $L$-orbits and the closed
one will be denoted by
$X_+$, which is naturally isomorphic to one component of the Grassmanian of
maximal isotropic vector subspaces in $W$.
Let  $\widehat X_+ \subset \S_+$ be the cone of $X_+$.

Recall \cite[p.102,103]{chevalley}
that there is an $L$-equivariant quadratic
map $\nu_2: \S_+ \rightarrow W$ with $\nu_2^{-1}(0) = \widehat X_+ $ and the
image is contained in the affine cone of $\cQ$.
This gives a rational map $\bar{\nu_2}: \p \S_+ \dasharrow \cQ.$

Let us denote by $CS_+$ the variety of ``complete 10-dimensional
spinors'' defined as the graph closure of the rational map
$\bar{\nu_2}$. This name comes from the cases $A,D$. By the
following lemma \ref{complete_spinor}, $CS_+ \to \p \S_+$ is
isomorphic to the blow-up  $\phi: Bl_{X_+}(\p \S_+) \to \p \S_+$
of $\p \S_+$ along $X_+$.

Let $Bl_1(\S_+) \to \S_+ $ be the blow-up of $\S_+$ along the
origin, then the exceptional fiber is $\p \S_+$ and $Bl_1(\S_+)$
is the total space of the tautological line bundle
$\mathcal{O}(-1)$ over $\p \S_+$. The strict transform of
$\widehat X_+$ is then the total space of the line bundle
$\mathcal{O}(-1)|_{X_+}$, which is smooth.  Claims (i) and (ii) of
the theorem follow immediately. Let $Bl_2(\S_+)$ be the blow-up of
$Bl_1(\S_+)$ along the strict transform of $\widehat X_+$. Then
the above discussions give that $Bl_2(\S_+)$ is just the total
space of the line bundle $\phi^* \0(-1)$ over $Bl_{X_+}(\S_+)
\simeq CS_+$. The second graph projection gives a morphism $CS_+
\to \cQ$, thus a morphism $\phi^* \0(-1) \to \cQ$. This is the
local model, and the globalization can be obtained as follows.

\smallskip

By the Bruhat decomposition theorem, all homogeneous spaces
(such as $X_P$) are covered by affine spaces, and
thus are locally isomorphic with their tangent spaces. In our case, we thus
have an isomorphism $U_P \simeq \S_+^*$ (where $U_P \subset X_P$ is a
suitable open subset). We deduce that that the restriction of $T^* X_P$ to
$U_P$ is isomorphic with $\S_+^* \times \S_+$.
Moreover, this isomorphism is $L$-equivariant.

By \cite{igusa}, there are three $L$-orbits in $\S_+$, which give via
this isomorphism the intersection of the three $G$-orbits in $T^* X_P$ with
$T^* U_P$.
Therefore, above $U_P$, the closure
of the $G$-orbits are trivial fibrations, and to understand the result of
the successive blow-ups, it is enough to understand this in one fiber.

Recall from  \cite{chaput} that $X_P$ has a minimal equivariant
projective embedding $X_P \subset \p V$, with $\dim V = 27$, and
that $X_Q$ naturally embeds in $\p V^*$. In terms of this
embedding, $W$ is a linear subspace of $V^*$ (in fact, it is the
orthogonal of the 17-dimensional affine tangent space
$\widehat{T_xX_P} \subset V$), and $\cQ = X_Q \cap \p W \subset \p
W$. Recall that $L$ is a Levi subgroup of $P$, so we can regard
$W$ as a $P$-representation. Now consider the vector bundle $G
\times^P W$ over $X_P$ (it is in fact the conormal bundle of $X_P$
in $\p V$) and the relative quadric over $X_P$ defined by $G
\times^P \cQ$. The previous construction gives a map $\tilde T_P
\rightarrow G \times^P \cQ$ over $X_P$. We have a natural
inclusion $G \times^P \cQ \to G \times^P (G/Q) \simeq G/P \times
G/Q$. The second projection of the latter yields our desired
morphism $\tilde T_P \rightarrow X_Q$. Its restriction to each
fiber, considered as a rational map, is an $L$-equivariant
rational map from a projectivised spinor representation of $L$ to
the 8-dimensional projective smooth quadric; by \cite[proposition
1.5]{chaput}, it coincides with $\pi_Q \circ \mu$.

By lemma \ref{lift} and \cite[theorem 4.1]{chaput}, this morphism
lifts to a morphism $\tilde T_P \rightarrow T_Q$. Assume that it
lifts further to a $G$-equivariant morphism $\tilde T_P
\rightarrow Bl_{k-1}(T_Q)$ for some $k \in \{1,2\}$, then one
checks that the preimage of $Bl_{k-1}(\overline{T_Q^{k-1}})$ is a
divisor, which implies that this morphism lifts to a
$G$-equivariant morphism $\tilde T_P \rightarrow Bl_k(T_Q)$. We
thus get a morphism $\tilde \mu : \tilde T_P \rightarrow \tilde
T_Q$. We have already checked that the restriction of $\tilde \mu$
to $\tilde T_P^r \simeq T_P^r$ to $\tilde T_Q^r \simeq T_Q^r$
identifies with the flop $\mu$.

We now consider the analogous morphism
$\tilde \mu_Q : \tilde T_Q \rightarrow \tilde T_P$. Since it also
coincides with the Mukai flop on the open orbit, it follows that
$\tilde \mu_Q \circ \tilde \mu : \tilde T_P \rightarrow \tilde T_P$ is the
identity on the open orbit; therefore it is the identity. So the
theorem is proved in type $E_6$.

\lpara

We finally explain how this proof adapts to the type $D_{2n+1}$.
Assume that $G$ is of type $D_{2n+1}$. Let $\C^{4n+2}$ be equipped
with a non-degenerate quadratic form. Let $\alpha \in G/P$ and
$L_\alpha$ the corresponding maximal isotropic subspace. As
before, the fiber $\tilde T_\alpha$ of $\tilde T_P \rightarrow
X_P$ over $\alpha$ is the result of the successive blowups of
$T_\alpha (G/P) \simeq \wedge^2 L_\alpha$ along the strict
transforms of the $P$-orbits in $T_\alpha (G/P)$. Recall that the
space of complete skew-symmetric forms in $L_\alpha^*$ is the
graph of the rational morphism
$$
\begin{array}{rcl}
\p (\wedge^2 L_\alpha) & \dasharrow & \p (\wedge^2 L_\alpha)
\times \p (\wedge^4 L_\alpha) \times \cdots \times
\p (\wedge^{2n} L_\alpha) \\
{} [\omega] & \mapsto & ([\omega],[\omega \wedge \omega],\ldots,
[\omega^{\wedge n}])
\end{array}
$$
By the following proposition \ref{complete_forms},
$\tilde T_\alpha$ identifies with
the total space of a line bundle over the space of complete
skew-symmetric forms in $L_\alpha^*$.

Now, a point in $\tilde T_\alpha$ therefore defines an element
$(x_1,\ldots,x_n) \in \p (\wedge^2 L_\alpha) \times \p (\wedge^2
L_\alpha) \times \cdots \times \p (\wedge^{2n} L_\alpha)$ in the
set of complete skew-symmetric forms in $L_\alpha^*$. The last
element $x_n$ belongs to $\p (\wedge^{2n} L_\alpha) = \p
L_\alpha^*$; thus it defines a hyperplane in $L_\alpha$. Since
there is a unique maximal isotropic subspace parameterized by an
element in $X_Q$ and meeting $L_\alpha$ along this hyperplane,
there is a natural morphism $\tilde T_\alpha \rightarrow X_Q$.
Then we can apply Lemma \ref{lift} and the previous arguments to
conclude. \qed

We now turn to the proof of the lemma used in the previous proof.
Recall that we defined
$$
CS_+ := \overline{ \{([s],[\nu_2(s)]) \mid  s \in \S_+ - \widehat X_+ \} }
\subset  \p \S_+ \times \cQ,
$$
where $\nu_2 : \S_+ \rightarrow W$ is a $Spin_{10}$-equivariant
quadratic map.

Recall that there exists a unique  symmetric bilinear map $\tilde
\nu_2: \S_+ \times \S_+ \rightarrow W$ such that $\tilde
\nu_2(s,s) = 2\nu_2(s)$. Recall that $W$ is the natural
10-dimensional representation of $Spin_{10}$; it is therefore
equipped with a scalar product that we denote by $\scal{ \cdot ,
\cdot}$.

\begin{lemm}
\label{equation_cs+} The variety $CS_+$ coincides with the variety
$$C:= \{([s],[x]) \in \p \S_+ \times \cQ | \forall t \in \S_+,
\scal{\tilde \nu_2(s,t),x} = 0 \}.$$
\end{lemm}
\begin{proo}
 For  $s \in  \S_+ - \widehat X_+$, let
$x=\nu_2(s) \in \widehat \cQ$. For any $t \in \S_+$, we have
$\nu_2(s+t \epsilon) \in \widehat \cQ$ for all $\epsilon$. Thus
$ \frac{d}{d \epsilon}|_{\epsilon = 0} \nu_2(s+t \epsilon) =
\tilde \nu_2(s,t) \in \widehat{T_x\cQ}=x^\bot,$ thus $CS_+ \subset
C$. On the other hand, the projection morphism $C \rightarrow \cQ$
is a locally trivial bundle with projective spaces as fibers, so $C$ is
irreducible. Note that if $([s],[x]) \in C$ with $\nu_2(s) \neq 0$, then
$[\nu_2(s)] = [x]$ because $\im(\tilde \nu_2(s, \cdot)) = x^\bot.$
This shows that over $\p \S_+ - X_+$, the two irreducible projective varieties
$CS_+$ and $C$ coincide. As $CS_+ \subset C$, we have the equality.
\end{proo}

\begin{lemm}
\label{complete_spinor}
$CS_+ \simeq Bl_{X_+}(\p \S_+)$.
\end{lemm}
\begin{proo}
By lemma \ref{equation_cs+}, $CS_+$ is the total space of a
locally trivial fibration with projective spaces as fibers, thus
$CS_+$ is a smooth variety. Moreover, the fiber over a point $[s]
\in X_+$ of the graph projection $CS_+ \to \p \S_+$ is then
isomorphic to $\p \im (\tilde \nu_2(s,.)) \simeq \p^4$. From this,
we see that the preimage of $X_+$ in $CS_+$ is a divisor. By the
universal property of blow-ups, it follows that there is a
morphism $f : CS_+ \rightarrow Bl_{X_+}(\p \S_+)$. Note that the
kernel of $\tilde \nu_2(s,.)$ is the affine tangent space
$\widehat {T_{[s]} X_+}$, and $f$ restricts to an isomorphism
between $\p \im (\tilde \nu_2(s,.))$ and $\p (\S_+/\widehat
{T_{[s]} X_+})$. Therefore, $f$ is a bijective birational
morphism. By Zariski's main theorem, it is an isomorphism.
\end{proo}

We now turn to the case of complete skew-forms. Let $n$ be an integer, and
let $V$ be a vector space of dimension $2n$ or $2n+1$.

\begin{defi}
\
\begin{itemize}
\item The variety of complete skew-forms on $V$ is the closure of
the graph (denoted by $\overline{\p \wedge^2 V^*}$) of the
rational map
$$\fonctionrat{\psi}{\p (\wedge^2 V^*)}{\p (\wedge^4 V^*) \times
\cdots \times \p (\wedge^{2n} V^*)}{[\omega]} {([\omega \wedge
\omega],\ldots,[\omega^{\wedge n}]).}$$  \item We denote by $Bl(\p
\wedge^2 V^*)$ the variety obtained by blowing up successively the
strict transforms of the different $SL(V)$-orbits in $\p \wedge^2
V^*$, from the smallest to the biggest.
\end{itemize}
\end{defi}

\begin{prop}
\label{complete_forms} There is an $SL(V)$-equivariant isomorphism
$\overline{\p \wedge^2 V^*} \simeq Bl(\p \wedge^2 V^*)$.
\end{prop}
\begin{proo}
This result on such a classical subject should be known to
specialists. We include a short proof here, taking advantage of
the recursive nature of complete skew forms, only because we were
not able to find it in the literature. In \cite[theorem
11.1]{tha}, a similar result is proved, but it is not precisely
what we need.

We first show by induction that $\overline{\p \wedge^2 V^*}$ and
$Bl(\p \wedge^2 V^*)$ are smooth, and that the pre-image of each
orbit in $\p \wedge^2 V^*$ is an irreducible divisor.  We choose a
subspace $L \subset V$ of dimension 2. Let $U$ denote the open
subset of $\p \wedge^2 V^*$ of elements $[\omega']$ such that
$\omega'_{|L} \not = 0$. Given $[\omega'] \in U$, let
$S_{[\omega']}$ be the orthogonal  of $L$ in $V$ with respect to
$\omega'$, which defines a morphism $p: U \to G(\dim V - 2,V).$
Let $q: \overline{\p \wedge^2 V^*} \to \p \wedge^2 V^* $ be the
natural projection and $g: q^{-1}(U) \to G(\dim V - 2,V)$ the
composition. Fix an element $S_{[\omega']} \in G(\dim V - 2,V),$
then $V = L \oplus S_{[\omega']}$.  By definition,
$g^{-1}(S_{[\omega']}) = \overline{\{([\omega], \cdots,
[\omega^{\wedge n}]) | \omega \in \wedge^2 V^*, \omega(L,
S_{[\omega']})=0  \}}.$

We fix a non-zero element $\lambda \in \wedge^2 L^*$.
Any element $[\omega]$ such that $ \omega(L, S_{[\omega']})=0$ can
be uniquely written as the class of $\lambda + \omega_0$, with
$\omega_0 \in \wedge^2 S_{[\omega']}^*$.
Note that $(\lambda + \omega_0)^{\wedge i} =
{\omega_0}^{\wedge i} + i \lambda \wedge
{\omega_0}^{\wedge(i-1)}.$ From this one can deduce that
$g^{-1}(S_{[\omega']})$ is
isomorphic to the total space of  a
line bundle $L$ over $\overline{\{([\omega_0], \cdots,
[\omega_0^{\wedge {n-1}}]) | \omega_0 \in \wedge^2 S_{[\omega']}^*\}},$
namely  the pull-back of the tautological line
bundle on $\p \wedge^2 S_{[\omega']}^*.$ Thus
$g: q^{-1}(U) \to G(\dim V - 2,V)$ is  a fibration with
fiber $L$. Thus the claim follows from induction. Similar
arguments prove the claim for $Bl(\p \wedge^2 V^*)$.

By the universal property of blow-ups, we obtain a surjective
morphism $f: \overline{\p \wedge^2 V^*}\to Bl(\p \wedge^2 V^*)$.
By Zariski's main theorem, $f$ has connected fibers, thus every
exceptional divisor of the morphism $\overline{\p \wedge^2 V^*}
\to \p \wedge^2 V^*$ is mapped birationally to an exceptional
divisor of $Bl(\p \wedge^2 V^*) \to \p \wedge^2 V^*$, which gives
that $Exc(f)$ has codimension $\geq 2$ if non-empty. As $f$ is a
birational map between two smooth varieties, thus $Exc(f)$ is
empty and $f$ is an isomorphism.
\end{proo}

We will now study in more detail the case of $G = E_6$, in
particular, we will give another more geometric proof of our
theorem in this case.  Let $Z$ denote the graph closure of the
flop $T_P \dasharrow T_Q$ and $R=P \cap Q$ the standard parabolic
subgroup corresponding to the roots $\alpha_1, \alpha_6$. We
denote by $\u(P), \u(Q)$ the nil-radical of the Lie algebras of
$P$ and $Q$ respectively.  Let $\n : = \u(P) \cap \u(Q)$. Then
$\n$ is an 8-dimensional $R$-representation and it is shown in
\cite[Theorem 6.1]{fu} that the graph closure $Z$ is isomorphic to
$G \times^R \n$. Furthermore the first graph projection $Z \to
T_P$ is isomorphic to the composition of natural maps $G \times^R
\n \simeq G\times^P (P \times^R \n) \to G \times^P (P \cdot \n)$.
In particular,  $P \cdot \n  = \u(P)$ is naturally identified to
$\S_+$.
\begin{lemm}\label{CS}
The natural map $\eta: P \times^R \p \n \to \p \S_+$ is the
blow-up $\phi: Bl_{X_+}(\p \S_+) \to \p \S_+$ of $\p \S_+$ along
$X_+$.  The natural map $G \times^R \p \n \to G \times^P
\p(\u(P))$ is the blow-up of $ G \times^P \p(\u(P))$ along the
unique $G$-closed orbit, say $ O_P$.
\end{lemm}
\begin{proo}
First we study the $R$-representation $\n$. Let $L$ be a Levi
factor of $R$. As one can see on the Dynkin diagram, $L$ contains
$Spin_8$, so let $M \simeq Spin_8$ be a subgroup of $L$. We
consider an element $x$ in the root space corresponding to the
highest root of $\e_6$. With the notations of \cite{bourbaki},
this highest root is $\tilde \alpha =
\alpha_1+2\alpha_3+3\alpha_4+2\alpha_5+\alpha_6+2\alpha_2$. If
$h_{\alpha_i}$ is the element of the Cartan subalgebra
corresponding to the simple root $\alpha_i$, it is therefore easy
to compute that $\tilde \alpha(h_{\alpha_i}) = 0,0,0,1$ for
$i=3,4,5,2$. It therefore follows that $x$, considered as an
element of the $M$-representation $\n$, has positive weight, and
is not a sum of two positive weights. Therefore, its weight is a
fundamental weight of $M \simeq Spin_8$, and since $\dim \n = 8$,
we deduce that $\n$ is either the natural or a spinor
representation of $M \simeq Spin_8$. It is thus a consequence of
the triality principle which exchanges these three representations
that there are three $R$-orbits in $\n$, of dimensions 0, 7 and 8
respectively.
This implies that in $\p \n$ there are only two
$R$-orbits.

Let $O$ be the closed $P$-orbit in $P \times^R \p \n
$, which is of codimension 1. The map $\eta$ contracts $O$ to
$X_+$, and it maps the open $P$-orbit isomorphically to the open
$P$-orbit in $\p \S_+$. By the universal property of blow-ups, we
have a morphism $\psi: P \times^R \p \n \to Bl_{X_+}(\p \S_+)$,
which induces a surjective (thus generically finite)  map $\psi_1:
O \to Exc(\phi)$. By Zariski's main theorem, the fibers of $\psi$
(thus also $\psi_1$) are connected, so $\psi_1$ is generically of
degree one, which implies that $\codim\ Exc(\psi) \geq 2$ if
$Exc(\psi) \neq \emptyset$. On the other hand, $\psi$ is a
birational morphism between two smooth varieties, so $Exc(\psi)$
is either empty or of pure codimension 1.  Thus $Exc(\psi)=
\emptyset$ and $\psi$ is an isomorphism. The second claim follows
immediately.

\end{proo}

\begin{lemm}
The birational map $G \times^P \p(\u(P)) \dasharrow G \times^Q
\p(\u(Q))$ is a family of $\p^5$-Mukai flops with center $ O_P$.
\end{lemm}
\begin{proo}
Note that $ O_P$ has two natural fibrations: $ O_P \to G/P$ and $
f:  O_P \to \p \0_{min} = G/P_2$, where $P_2$ is the maximal
parabolic associated to the root $\alpha_2$.
Since moreover $O_P$ is closed, this shows that
$O_P$ is the incidence variety $G/(P \cap P_2)$. Thus $f$ is a
locally trivial $P_2/(P \cap P_2) \simeq \p^5$-bundle. Let $\tilde
O_P$ be its pre-image in $G \times^P (\u(P)-0)$, then $\tilde f:
\tilde O_P \to \0_{min}$ is again a $\p^5$-bundle, thus the normal
bundle of $\tilde O_P$ is isomorphic to $\Omega_{\tilde f}$ (see
for example \cite{huy}, section 3). Now it follows that the normal
bundle of $ O_P$ in $G \times^P \p(\u(P))$ is isomorphic to
$\Omega_{f}$. By Lemma \ref{CS}, the blowing-up of $G \times^P
\p(\u(P))$ along $\bar O_P$ gives $G \times^R \p \n$, while the
latter is also the blowing-up of $G \times^Q \p(\u(Q))$ along
$O_Q$, concluding the proof.
\end{proo}

We now give a more geometric proof of theorem \ref{resolution} in
the case $E_{6,I}$~:

\begin{theo}\label{E6}
Consider the stratified Mukai flop of type $E_{6, I}, \mu: T_P \dasharrow T_Q.$
 Let $Bl_1(T_P)\to T_P$ (resp. $Bl_1(T_Q) \to T_Q$) be the blow-up
of $T_P$ (resp. $T_Q$) along the zero section and $Bl_1(\overline
T_P^1)$ (resp. $Bl_1(\overline T_Q^1)$) the strict transform of
the closure of $\tilde O_P$ (resp. $\tilde O_Q$). Then the
birational map $Bl_1(T_P) \dasharrow Bl_1(T_Q)$ is a family of
$\p^5$-Mukai flops with center $Bl_1(\overline T_P^1)$. In
particular, there exists a $G$-equivariant isomorphism $\tilde
\mu: \tilde T_P \simeq \tilde T_Q$ which identifies with $\mu$
over the open $G$-orbits.
\end{theo}
\begin{proo}
Let $\bar \0$ be the nilpotent orbit closure which is resolved by
$T_P$ and $T_Q$.  Let $Bl_1(\bar \0)$ be the blow-up of $\bar \0$
along the origin. By the universal property of blow-ups, we obtain
morphisms $Bl_1(T_P) \to Bl_1(\bar \0) \leftarrow Bl_1(T_Q). $ To
simplify the notations, we let $Y_P=G \times^P \p(\u(P))$ and
$Y_Q= G \times^Q \p(\u(Q))$. Then $Bl_1(T_P)$ (resp. $Bl_1(T_Q)$,
$Bl_1(\bar \0)$) is just the total space of the tautological line
bundle $\0_{Y_P}(-1)$ (resp. $\0_{Y_Q}(-1)$,  $\0_{\p \bar
\0}(-1)$). The birational map $Bl_1(T_P) \dasharrow  Bl_1(T_Q)$ is
just the pull-back via the morphism $\0_{\p \bar \0}(-1) \to \p
\bar \0$ of the birational map $Y_P \dasharrow Y_Q$. Now the claim
follows from the precedent lemma.

\end{proo}

Now we will relate the variety $\tilde T_P$ to the graph closure $Z$ of the flop $\mu$.
 Let $S \simeq G/R$ denote the closed $G$-orbit in $Z$.
By the precedent theorem,  $\tilde T_P$ gives a resolution of the flop $\mu: T_P \dasharrow T_Q$, thus there exists
 a birational morphism $\gamma: \tilde T_P \to Z$.
\begin{coro}\label{graphe}
The morphism $\gamma: \tilde T_P \to Z$ is isomorphic to the
blow-up of $Z$ along $S$.
\end{coro}
\begin{proo}
Note that $Z$ is smooth, so $Exc(\gamma)$ is of pure codimension
1. The complementary of the open $G$-orbit in $\tilde T_P$
consists of two irreducible divisors: $\tilde T_P^0$ and $\tilde
T_P^1$. By Lemma \ref{CS}, $\gamma$ maps $\tilde T_P^1 -  \tilde
T_P^0$ isomorphically to its image, which is the blow-up of the
27-dimensional $G$-orbit in $T_P$. This implies that $supp(
Exc(\gamma)) = \tilde T_P^0$.
 Now a similar argument as that in the proof of Lemma \ref{CS} proves the corollary.

\end{proo}

To summarize, the geometric picture of the  flop $\mu: T_P
\dasharrow T_Q$ is  as follows: the 27-dimensional $G$-orbit is a
$\p^5$-bundle over the 22-dimensional orbit in  the nilpotent
orbit  closure $\overline{\0}$. The flop  $\mu$ when restricted to
the complementary of the zero section in $T_P$ is just a standard
Mukai flop of  $\p^5$-bundles. This restricted flop is resolved by
just one blow-up, which is given by $G \times^R (\n-\{0\})$. Thus
$\mu$ can be regarded as a degeneration of a 22-dimensional family
of $\p^5$-Mukai flops. Note that such a picture is very similar to
stratified Mukai flops: $T^* Gr(2, V) \dasharrow T^* Gr(2, V^*).$

\section{Deformations of stratified Mukai flops}

In this section,  we will consider only the $E_{6, I}$ case, which
share many analogous properties with the standard Mukai flop (see
for example  \cite{chaput}). As before, let $P:=P_1$ and $P_6$ be
the maximal parabolic subgroups of $E_6$. They have conjugate Levi
subgroups, so we can choose $Q$  in the conjugacy class of $P_6$
such that $P$ and $Q$ have the same Levi subgroup $L$. Note that
the center $\mathfrak{c}$ of the Lie algebra of $L$ is
one-dimensional.

\begin{lemm}
\label{stabilisateur}
$\mathfrak{c} \oplus \u(P)$ is stable under $P$, and the stabilizer of the
line $\mathfrak{c}$, considered as an element in $\p \g$, is $L$.
\end{lemm}
\begin{proo}
Let $U$ denote the unipotent radical of $P$. According to the Levi
decomposition theorem, we have $P=LU$. Therefore, if $p=lu$, and
$x \in \mathfrak{c}$, we have $p.x \equiv l.x \equiv x$ modulo
$\u(P)$, proving the first point.

For the second point, denote $\mathfrak{z}(\mathfrak{c})$ the set
of elements $z \in \g$ such that $\forall c \in
\mathfrak{c},[z,c]=0$ and $Stab(\mathfrak{c}) = \{ g \in G :
g.\mathfrak{c} \subset \mathfrak{c} \}$. Note that $\mathfrak{l}
\subset \mathfrak{z}(\mathfrak{c})$, and for reason of equal
dimension,  we  have the equality . Therefore,
$Stab(\mathfrak{c}) \subset Stab(\mathfrak{z}(\mathfrak{c}))
=Stab(\mathfrak{l})$. Now let $\plie$ and $\q$ be the two
parabolic subalgebras containing $\mathfrak{l}$, and let $P,Q$ be
the corresponding parabolic subgroups. We have $Stab(\mathfrak{l})
\subset Stab(\plie) \cap Stab(\q) = P \cap Q = L$.
\end{proo}

Let $T_P = T^*(G/P)$ and $T_Q = T^*(G/Q)$ be the
cotangent spaces. We have the one-dimensional smooth flat
deformations of these two varieties: $E_P: = G \times^P
(\mathfrak{c} + \u(P))$ and $E_Q: = G \times^Q (\mathfrak{c} +
\u(Q))$.

\begin{lemm}
There exists a $G$-equivariant
birational map $\psi: E_P \dasharrow E_Q$ which deforms the
Mukai flop $\mu: T_P \dasharrow T_Q$ to isomorphisms outside the zero fibers.
\end{lemm}
\begin{proo}
Let $z \in \mathfrak{c}$ be a non-zero element. By the previous lemma, its
stabilizer in $G$ is $L$ itself. Let $U$ be the unipotent part of
$P$, then $P \cdot z = U \cdot z$.  Note that $U \cdot z$ is
closed and has  $ z + \u(P)$  as the tangent space at $z$, which
gives $z + \u(P) = P \cdot z,$ thus $G\cdot (z+\u(P))= G \cdot z$
and $\cN:=G \cdot (\mathfrak{c}+\u(P)) =  G \cdot
(\mathfrak{c}+\u(Q)).$

We have a natural $G$-equivariant  projective morphism $e_P: E_P
\to \cN$. Lemma \ref{stabilisateur} implies that $e_P$ is an
isomorphism outside the singular locus of $\overline{\0}:=G \cdot
\u(P) = G \cdot \u(Q)$, so $e_P$ is a small resolution. In a
similar way, we have a small resolution $e_Q: E_Q \to \cN$. Now
the claim follows.
\end{proo}

Our next task is to resolve this birational map. Recall we have a
stratification $T_P^0 \subset \overline{T_P^1} \subset T_P$. We
will regard these varieties as subvarieties in $E_P$. Let
$Bl_1(E_P)$ be the blow-up of $E_P$ along $T_P^0$ and  $\tilde
E_P$ the blow-up of $Bl_1(E_P)$ along the strict transform of $
\overline{T_P^1}$. Similarly we can get $\tilde E_Q$. The next
theorem is analogous to Theorem \ref{resolution}.
\begin{theo} \label{def-reso}
There exists a unique $G$-equivariant isomorphism $\tilde E_P \simeq \tilde E_Q$.
\end{theo}
\begin{proo}
Let $\0_1$ be the 22-dimensional nilpotent orbit in
$\overline{\0}$. The pre-image of $\0_1$ under the symplectic
resolution $\Pi_P: T_P \to \overline{\0}$ is a $\p^5$-bundle, i.e.
$\Pi_P^{-1} (\0_1) = \p(F)$ for some vector bundle $\phi: F \to
\0_1$ of  rank 6. By our discussions in the precedent section, the
variety $\Pi_Q^{-1}(\0_1)$  is isomorphic to $\p(F^*)$.

Let $v \in H^2(T^*G/P) \simeq H^2(G/P)$
be the Kodaira-Spencer class of the deformation $E_P \to
\mathfrak{c}.$ Note that $v \neq 0$ since the deformation is
non-trivial. The Picard group of $G/P$ has rank one, thus  either
$v$ or $-v$ is  $\pi_P$-ample. In both cases, the restriction of
$v$ to a fiber of $\phi$ is non-trivial. By the proof of
\cite[Lemma 3.6]{huy}, the normal bundle $\cN_{\Pi_P^{-1}(\0_1)
\mid T_P}$ is isomorphic to $\phi^*(F^*) \otimes \0_\phi(-1).$ Let
us denote by $p_P: Bl_1(\p E_P) \to \p E_P$ the blow-up of $\p
E_P$ along $\p T_P^1$.  Then we have an isomorphism $Bl_1(\p E_P)
\simeq Bl_1(\p E_Q).$ In particular, the birational map: $\p E_P
\dasharrow \p E_Q$ is a family of $\p^5$-standard flops. Now the
claim follows from the same arguments of the proof of Theorem
\ref{E6}.

\end{proo}

\section{Chow groups}
We will denote by $\tilde E$ the variety $\tilde E_P \simeq \tilde
E_Q$.  Let $f^P: \tilde E \to E_P$ and $f^Q: \tilde E \to E_Q$ be
the two natural morphisms. Let $\tilde T$ be the central fiber of
the family $\tilde E \to \C$. Note that $\tilde T$ has three
irreducible components, one of which is $\tilde T_P$. Let  $f_0^P:
\tilde T \to T_P$ and $f_0^Q: \tilde T \to T_Q$ be the natural
morphisms.  We define two homomorphisms $\Psi=f^Q_* \circ (f^P)^*:
CH(E_P) \to CH(E_Q)$ and $\Psi_0=(f_0^Q)_* \circ (f_0^P)^*: CH(T_P)
\to CH(T_Q)$. Here $CH(X)$ denotes the Chow ring with integer
coefficients of the variety $X$.
\begin{theo}
The maps $\Psi, \Psi_0$ are both isomorphisms.
\end{theo}
\begin{proo}
For the proof of the first part, we will follow the proof of Theorem
2.1 in \cite{LLW}, where a similar result is proved for standard
flops.

Let $W_P$ be any $k$-dimensional sub-variety in $E_P$, by Chow's
moving Lemma (see for example \cite{Ful}, section 11.4), up to
replacing $W_P$ by an equivalent cycle, we can assume that $W_P$
intersects properly with the cycle $\bar T_P^1 + T_P^0$. Applying
twice the blow-up formula (see \cite{Ful}, Theorem 6.7), we obtain
that $(f^P)^*(W_P) = \tilde W$, where $\tilde W$ is the strict
transform of $W_P$ under the birational map $f^P$, which gives
that $\Psi(W_P) = W_Q$, where $W_Q$ is the strict transform of
$W_P$ under the birational rational map $\psi: E_P \dasharrow
E_Q$. Remark that $W_Q$ intersects  no longer properly with the
cycle $\bar T_Q^1 + T_Q^0$.

Note  that for any irreducible component $C$ of  $ W_Q \cap \bar
T_Q^1 $, there exists an irreducible component $B$ of the
intersection $W_P \cap \bar T_P^1$ such that $C \subset
e_Q^{-1}(e_P(B))$. Let $(f^Q)^* W_Q = \tilde W + \Sigma F_C $,
where $F_C \subset (f^Q)^{-1}e_Q^{-1}e_P(B)$. Note that we have
$\dim B = \dim W_P - \codim \bar T_P^1. $
 For a generic point $s \in e_P(f^P(F_C))$, we
have
$$ \dim  F_{C,s} = \dim F_C - \dim e_P(B)  \geq  \dim W_P - \dim B = \codim \bar T_P^1 =6.$$
However a generic fiber of $e_P |_{\bar T_P^1}$ is 5-dimensional,
thus $F_{C,s}$ contains positive dimensional fibers of $f^P$, thus
$f^P_* (F_C) =0$, which gives that $f^P_* (f^Q)^*$ is an inverse
of $\Psi$, thus the claim of  the first part.

Let $i_P: T_P \to E_P, i_Q: T_Q \to E_Q $ and $i: \tilde T \to
\tilde E$ be the natural inclusions. The notations are summarized
in the following diagram:
\begin{diagram}[size=0.8cm]
    &               &   T_P   & \rTo^{i_P}  & E_P  &  &  \\
    & \ruTo^{f_0^P} &   &   \ruTo^{f^P} & &\rdTo^{e_P}  & \\
\tilde{T} & \rTo^i  & \tilde{E}  &  & & & \cN \\
 & \rdTo^{f_0^Q} &   &   \rdTo^{f^Q} & &\ruTo^{e_Q}  & \\
&               &   T_Q   & \rTo^{i_Q}  & E_Q  &  &  \\
\end{diagram}

 We have the following
diagrams:
\[
\begin{CD}
CH(E_P) @>(f^P)^*>> CH(\tilde E) @>(f^Q)_*>> CH(E_Q) \\
@Vi_P^*VV  @Vi^*VV  @VVi_Q^*V  \\
CH(T_P) @>(f^P_0)^*>> CH(\tilde T) @>(f^Q_0)_*>> CH(T_Q)
\end{CD}
\]
The first diagram is commutative since $f^P \circ i = i_P \circ
f^P_0: \tilde T \to E_P.$  Now we show that the second diagram is
also commutative. In fact, for any $[Y] \in CH(\tilde E),$ we may
assume it intersects properly with $\tilde T$ by using Chow's
moving lemma. Note that $f^Q$ is an isomorphism outside $\tilde
T$, thus $(f^Q)_* [Y] = [f^Q(Y)]$. It follows that $i_Q^* (f^Q)_*
[Y] = [f^Q(Y) \cap T_Q].$ Similarly we have $(f_0^Q)_* i^* [Y] =
[f_0^Q(\tilde T \cap Y)]$, which gives $(f_0^Q)_* i^*=i_Q^*
(f^Q)_*$. It follows  that $i_Q^* \circ \Psi = \Psi_0 \circ
i_P^*$. Now the claim follows from the fact that both $i_Q^*$ and
$i_P^*$ are both isomorphisms.
\end{proo}

\begin{rem}
(i). As shown in \cite{LLW}, this also proves that $\Psi$
(resp.  $\Psi_0$) induces an equivalence between the motives $[E_P]$
and $[E_Q]$ (resp. $[T_P]$ and $[T_Q]$).

(ii). The same proof also works for stratified Mukai flops on
Grassmanians $T^* Gr(2, V) \dasharrow T^* Gr(2, V^*),$ which has
been proved using a different method  by Namikawa (\cite{Nam}) on the level of
$K$-groups. His proof morally works for stratified Mukai flops of type $D$.
So the only unsolved case is the stratified Mukai flops of type $E_{6, II}$.
\end{rem}

\begin{exam}\label{ex1}
To simplify the notations, we denote by $b_P: B_P \to T_P$ the
blow-up along $X_P$ and $F_P$ its exceptional divisor.
Let
$\Sigma_P \subset T_P$ be the closure of the 27-dimensional
$G$-orbit and $\tilde \Sigma_P$ its proper transform under $b_P$.
Let $h_P: \tilde T_P \to B_P$ the blow-up along $\tilde
\Sigma_P$. Similar notations will be used on the dual-side.
Finally
let $T_P \xleftarrow{g_P} \tilde T_P \simeq \tilde T_Q
\xrightarrow{g_Q} T_Q$ be the compositions of blow-ups.
Now we will show that the naturally
defined morphism $\Phi:=(g_Q)_* \circ (g_P)^*: CH(T_P) \to
CH(T_Q)$ is not an isomorphism. Note that this situation is very
similar to that of usual Mukai flops (see for example
\cite[section 2]{Nam}, \cite[section 6]{LLW}).

By the blow-up formula, we have $(b_P)^*(\Sigma_P) = \tilde
\Sigma_P + \Gamma_P$, where $\Gamma_P$ is a cycle supported on
$F_P$. Using Chow's moving lemma (with ambient space $F_P$) , we
can assume that $\Gamma_P$ intersects properly with $F_P \cap
\tilde \Sigma_P.$ Then $\gamma:=(h_Q)_* \circ (h_P)^* \Gamma_P$
will be a cycle supported on $F_Q$. Note that $b_Q(F_Q)$ is
16-dimensional, so $(b_Q)_* \gamma =0$. Note that the birational
map $B_P \dasharrow B_Q$ is a family of $\mathbb{P}^5$ Mukai
flops. One deduces that (see Example 6.5 \cite{LLW}) that $(h_Q)_*
\circ (g_P)^* (\Sigma_P) = 5  \tilde \Sigma_Q + \gamma$, which
gives that $\Phi(\Sigma_P) = 5 \Sigma_Q$. Similarly for
$\Phi':=(g_P)_* \circ (g_Q)^*: CH(T_Q) \to CH(T_P)$, we get
$\Phi'(\Sigma_Q) = 5 \Sigma_P$. Thus $\Phi$ and $\Phi'$ are not
isomorphisms if $\Sigma_P$ is non-zero in $CH^5(T_P).$ This can be
seen as follows.

 Let $q: \p(T^*X_P \oplus \0) \to X_P$ be the natural projection
and $\bar \Sigma_P$ the closure of $\Sigma_P$. By \cite[Prop.
3.3]{Ful}, we have $(\pi_P^*)^{-1}(\Sigma_P) = q_*(c_{16}(\xi)
\cap \bar \Sigma_P),$ where $\xi$ is the universal quotient bundle
(of rank 16) of $q^*(T^*(X_P) \oplus \0).$ Note that $ c_{16}(\xi)
= \sum_{i=0}^{16} c_1(\0(1))^i \cap c_{16-i}(q^* T^*X_P).$ Thus we
get $(\pi_P^*)^{-1}(\Sigma_P) =\sum_{i=0}^{16} q_*( c_{16-i}(q^*
T^*X_P) \cap c_1(\0(1))^i \cap \bar \Sigma_P) = \sum_{i=0}^{16}
c_{16-i}(T^*X_P) \cap q_*(c_1(\0(1))^i \cap \bar \Sigma_P).$ For
reason of dimension, $q_*(c_1(\0(1))^i \cap \bar \Sigma_P)=0$ if
$i \neq 5$, and $q_*(c_1(\0(1))^5 \cap \bar \Sigma_P)= d\ [X/P],$
where $d$ is the degree of the closure of $\hat X_+$.
 From this, we deduce
that $(\pi_P^*)^{-1}(\Sigma_P) =d\ c_5(T^*X_P)$, which is non-zero
since $c_5(T^*X_P) \neq 0$ (this can be deduced for example from
the calculus made in \cite[section 7]{IL}).
\end{exam}

\begin{exam}\label{ex2}
Let $Z$ be the graph of the flop $T_P \dasharrow T_Q$ and
$T_P \xleftarrow{q_1} Z \xrightarrow{q_2} T_Q$ the two graph projections.
Let $\gamma: \tilde T_P \to Z$ be the natural morphism, which is a blow-up
along a smooth subvariety (Corollary \ref{graphe}).
 By proposition 6.7 (b) (\cite{Ful}),
the map $\gamma_* \circ \gamma^*: CH(Z) \to CH(Z)$ is the identity.

Now consider the morphism $\Gamma:=(q_2)_* \circ q_1^*: CH(T_P) \to CH(T_Q).$
Then $\Gamma =(q_2)_* \circ \gamma_* \circ \gamma^* \circ q_1^* = \Phi$.
Thus the morphism $\Gamma$ induced by the graph of the flop does not give an
isomorphism.
\end{exam}

\quad \\
Labo. J. Leray, Facult\'e des sciences, Universit\'e de NANTES \\
             2, Rue de la Houssini\`ere,  BP 92208
             F-44322 Nantes Cedex 03 - France \\
             chaput@math.univ-nantes.fr \\
             fu@math.univ-nantes.fr


\begin{thebibliography}{99}

\bibitem[Bou]{bourbaki}    N. Bourbaki,
                        {\it  Groupes et alg\`ebres de Lie, chapitre 6},
                        Hermann, 1968.

\bibitem[Cha]{chaput} P.-E. Chaput,
               {\it  On Mukai flops for Scorza varieties}.
               math.AG/0601734.

\bibitem[Che]{chevalley}   C. Chevalley,
                    {\it  The algebraic theory of spinors and Clifford algebras}.
                    Collected works. Vol. 2,
                    Springer-Verlag, Berlin, 1997.


\bibitem[C-P]{deconcini} C. De Concini, C. Procesi.   {\it  Complete symmetric
varieties}, in  Invariant theory (Montecatini, 1982),  1--44,
Lecture Notes in Math., {\bf 996}, Springer, Berlin, 1983.

\bibitem[Fau]{faulkner} J.R. Faulkner,  {\it Octonion planes defined by quadratic Jordan algebras}.
Memoirs of the American Mathematical Society, No. {\bf 104}, 1970.

\bibitem[Fu]{fu} B. Fu,
            {\it  Extremal contractions, stratified Mukai flops and Springer maps}.
            math.AG/0605431.
\bibitem[Ful]{Ful}
W. Fulton,  {\it Intersection theory}. Ergebnisse der Mathematik
und ihrer Grenzgebiete (3)  2. Springer-Verlag, Berlin, 1984.
\bibitem[Huy]{huy}
D. Huybrechts,  {\it  Birational symplectic manifolds and their
deformations}.
  J. Differential Geom.  45  (1997),  no. 3, 488--513.
\bibitem[Igu]{igusa} J. Igusa,
 {\it A classification of spinors up to dimension twelve}, Amer. J. Math.
{\bf 92} (1970) 997--1028.
\bibitem[I-L]{IL}
A. Iliev; L. Manivel,  {\it The Chow ring of the Cayley plane},
Compos. Math.  141  (2005),  no. 1, 146--160.
\bibitem[L-L-W]{LLW} Y.-P. Lee, H.-W. Lin, C.-L. Wang,
  {\it  Flops, motives and invariance of quantum rings}, math.AG/0608370.
\bibitem[Mar]{markman} E. Markman,
             {\it  Brill-Noether duality for moduli spaces of sheaves on $K3$ surfaces}.
         J. Algebraic Geom.  10  (2001),  no. 4, 623--694.
\bibitem[Na1]{Nam} Y. Namikawa,
 {\it Mukai flops and derived categories. II}, in {\em  Algebraic
structures and moduli spaces},   149--175, CRM Proc. Lecture
Notes, 38, Amer. Math. Soc., Providence, RI, 2004.
\bibitem[Na2]{Na2}
Y. Namikawa,  {\it  Birational geometry of symplectic resolutions
of nilpotent orbits}, to appear in {\em Advanced Studies in Pure
Mathematics}, see also math.AG/0404072 and math.AG/0408274.

\bibitem[Tha]{tha} M.  Thaddeus,
            {\it Complete collineations revisited.}
        Math. Ann. {\bf 315}, 469--495 (1999).

\end{thebibliography}
\end{document}